\documentclass[11pt]{article}
\usepackage{amsfonts}
\usepackage{latexsym}   
\usepackage{amsmath,amssymb,amsthm,latexsym}

\input cyracc.def
\begin{document}
\newtheorem{theorem}{Theorem}[section] \newtheorem{lemma}[theorem]{Lemma}
\newtheorem{definition}[theorem]{Definition}
\newtheorem{example}[theorem]{Example}
\newtheorem{remark}[theorem]{Remark}
\newtheorem{corollary}[theorem]{Corollary}
\newtheorem{proposition}[theorem]{Proposition}
\newtheorem{problem}[theorem]{Problem}

\baselineskip 0.25in

\title{\Large\bf   Atomic Effect Algebras with the Riesz Decomposition
Property}
\author{{Anatolij Dvure\v{c}enskij$^{1,2}$, Yongjian Xie$^{1,3}$\thanks{E-mail: yjxie@snnu.edu.cn},}\  {}\\
 {\small $^1$Mathematical Institute, Slovak Academy of Sciences, \v{S}tef\'{a}nikova 49, SK-814 73 Bratislava, Slovakia }\\
 {\small $^2$Depar.
Algebra  Geom.,  Palack\'{y} Univer.,
CZ-771 46 Olomouc, Czech Republic,  dvurecen@mat.savba.sk}\\
 {\small $^3$College of Mathematics and Information  Science, Shaanxi Normal University, Xi'an, 710062, China}\\
 }

\date{}
\maketitle
\begin{center}
\begin{minipage}{140mm}

\begin{abstract}
We discuss  the relationships
between  effect algebras with the Riesz Decomposition Property
and partially ordered groups with interpolation. We show that any
$\sigma$-orthocomplete atomic effect algebra with the Riesz
Decomposition Property is an MV-effect algebras, and we apply this result for pseudo-effect  algebras and for states.
\end{abstract}



\end{minipage}
\end{center}

\section{Introduction and basic definitions}

Effect algebras were introduced by Foulis and Bennett \cite{FouBen94} for the study
of logical foundations of quantum mechanics.
Independently, K\^{o}pka and Chovanec \cite{KopCho94} introduced essentially
equivalent structures called D-posets. If in
an effect algebra $(E;+,0,1),$ we  define a partial binary
difference operation $-$ as follows: for $a,b\in E$, $a-
b=c$ if and only if $b+ c$ exists in $E$ and $a=b+ c,$
then the algebraic system $ (E;-,0,1)$ is a D-poset
\cite{DvuPul00}. Effect algebras are a common generalization of
several well-established algebraic structures, in particular of
orthomodular lattices, orthomodular posets, orthoalgebras   and MV-algebras \cite{DvuPul00}.

The most important example of effect algebras is the system $\mathcal E(H)$ of all Hermitian operators of a (real, complex or quaternionic)  Hilbert space $H$ that are among the zero and the identity operator. $\mathcal E(H)$ is used for modeling unsharp observables via POV-measures in measurements in  quantum mechanics.

In 1958, Chang \cite{Cha} introduced  MV-algebras to prove the
completion of the \L ukasiewicz  propositional  logic \cite{Cig00}.
MV-algebras play an important role in many fields of mathematics
\cite{Cig00,DvuPul00}. Especially, MV-algebras have appeared in
effect algebras in many ways: Mundici showed that starting from any
AF C$^{*}$-algebra we can obtain a countable MV-algebra, and
conversely, any countable MV-algebra can be derived in such a way
\cite{Cig00,DvuPul00}. Ravindran \cite{Rav98} proved that $\Phi$-symmetric
effect algebras  are exactly MV-algebras, and also Boolean D-posets
of Chovanec and K\^{o}pka are MV-algebras \cite{DvuPul00}.
Especially,  Rie\v{c}anov\'{a} \cite{Rie00} has  proved that
every lattice-ordered effect algebra  $(E;+,0,1)$ is an
MV-algebra $(E; \oplus, ^{\prime},0, 1)$ if every pair of elements
of the effect algebra $E$ is compatible. Indeed, if we define a binary addition operation $\oplus$ on $E$ as
follows, for any $a,b\in E,$ $a\oplus b:=a+(a^{\prime}\wedge b)$, and
an unary operation $^{\prime}$ as follows: for any $a\in E$,
$a^{\prime}:=1- a.$  Conversely, in any MV-algebra $(E;
\oplus,$ $^{\prime},$ $0, 1)$, if we define a partial binary addition
operation $+$ on $E$ as follows: $a+ b$ exists if and only
if $a\leqslant b^{\prime}$ and in such a case $a+ b=a\oplus b,$ then the
algebraic system $(E;+,0,1)$ is a lattice-ordered effect algebra. We recall that any MV-algebra is also called an  MV-effect algebra, \cite{Jen04}.

Effect algebras with the Riesz decomposition property (RDP) form an important class of effect algebras. An effect algebra with RDP is
always an interval in an Abelian partially ordered group \cite{DvuPul00}. Every MV-effect algebra satisfies RDP. On the other hand, effect algebras with RDP are not necessarily MV-effect algebras.
However, every finite effect algebra with the Riesz decomposition property is   an MV-effect algebra \cite{BeFo95}. In this paper, we will continue in the study of the conditions  when effect algebras with RDP are MV-effect algebras.

The paper is organized as follows. In Section 2, we review some basic definitions and facts on effect algebras. In Section 3,  relationships between effect algebras with RDP and  partially
ordered Abelian groups with interpolation are discussed. In Section 4, we prove that any $\sigma$-orthocomplete atomic effect algebra with  RDP is also an MV-effect algebra.  Finally, in Section 5, we apply the results from Section 4 to a noncommutative generalization of effect algebras, called pseudo-effect  algebras, to show when they are effect algebras and to describe the state space of such effect algebras.

\section  {Basic definitions and facts}

\begin{definition}\label{def:effectalgebra}
{\rm \cite{FouBen94} An    {\it effect algebra}  is a system $(E;+,0,1)$ consisting of a set $E$ with two special elements $0$ and $1$, called the zero and the unit, and with  a partially binary operation $+$ satisfying the following conditions for
all $a,b,c\in E:$

\begin{itemize}
\item[{(E1)}] If $a+ b$ is defined, then $b+ a$ is defined and $a+ b=b+ a.$

\vspace{-2mm}\item[{(E2)}] If $ a+ b$ is defined and $(a+  b)+ c $ is defined, then    $ b+ c$ and $a+ ( b+ c)$ are defined, and   $(a+ b)+ c =a+ ( b+ c).$

\vspace{-2mm}\item[{(E3)}] For any $a\in E,$  there exists a unique $b\in E $ such that $a+ b$ is defined and $a+ b=1.$

\vspace{-2mm}\item[{(E4)}] If $ a+ 1$ is defined, then $a=0.$
\end{itemize}
}
\end{definition}

Let  $a$ be an element of an effect algebra $E$ and $n\geqslant0$ be an integer.
We define $na=0$ if $n=0,$ $1a=a$ if $n=1,$ and $na=(n-1)a+ a$ if $(n-1)a$ and $(n-1)a+ a$ are defined in $E$.  We define the {\it isotropic index} $\imath(a)$ of the element $a$, as the maximal nonnegative number $n$ such that $na$ exists. If $na$ exists for
every  integer $n$, we say that $\imath(a)=+\infty.$

\begin{remark}\label{re:bioper}
{\rm \cite{DvuPul00} Let $(E; +, 0, 1)$ be   an effect algebra.

(i)\ Define a partial binary relation $\leqslant$ on $E$ by $a\leqslant b$
if,  for some $c\in E,$ we have $c+ a=b.$ Then $(E; \leqslant, 0, 1)$ is
a poset, and $0\leqslant a\leqslant 1$ for each $a \in E.$ Furthermore, if $(E; \leqslant, 0, 1)$ is a lattice, then
we say that  $(E; +, 0, 1)$ is a lattice effect algebra.

(ii)\ Define a binary relation $\bot$ on $E$ by $ a\bot b$ if and
only if $a+ b$ exists in $E.$

(iii)\ Define a partial binary operation $-$ on $E$ by
$c- b=a$ if and only if $a+ b=c.$ Then the algebraic
system $(E; -, 0, 1)$ is a D-poset, \cite{DvuPul00}.
}
\end{remark}

For a comprehensive review on effect algebras, see \cite{DvuPul00}, where also unexplained notions from this paper can be found.

Let $(E;\leqslant)$ be a poset and let $a,b\in E$ be two elements such that $a \leqslant b.$  Then we define an interval $E[a,b]:=\{c \in E \mid a\leqslant c \leqslant b\}.$

We recall that a group $(G;+,0)$ written additively is a {\it partially ordered group} (po-group for short) if $a \leqslant b$ implies $c+a+d \leqslant c+b+d$ for all $c,d \in G.$  If $G$ with respect to $\leqslant$ is a lattice, we call $G$ a lattice-ordered group ($\ell$-group for short).

We denote by $G^+:=\{ g\in G \mid 0 \leqslant g\}$ the positive cone of $G.$ A po-group is {\it directed} if, for any $g_1,g_2\in G,$ there is an element $h\in G$ such that $g_1,g_2 \le h.$

If $G$ is a po-group and $u \in G^+,$ then the interval $G^+[0,u]$ can be converted into an effect algebra if we say that, for $a,b \in G^+[0,u],$ $a+b$ is defined in $G^+[0,u]$ iff the group addition $a+b$ is  in $G^+[0,u]$ and our addition $a+b$ coincides then with the group addition. Then $(G^+[0,u]; +,0,u)$ is an effect algebra.  Every effect algebra $E$ which is isomorphic with some $G^+[0,u],$ where $G$ is a po-group with strong unit $u$, is said to be an {\it interval effect algebra}.

\begin{definition}\label{de:atom}
{\rm  (i) An element $a$ of a poset $E$  with the least element $0$ is called an {\it atom}, if the interval $E[0,a]=\{x\in E\mid 0\leqslant x\leqslant a\}$ equals  the set $\{0,a\}$.

(ii) An   effect algebra $E$ is called   {\it atomic} if, for any nonzero $x$ of $E$, there exists an atom $a$  in $E$ such that $a\leqslant x$.
 }
\end{definition}

\begin{definition}\label{de:t-norm}
{ \rm \cite{DvuPul00} An effect algebra $(E;+,0,1)$ has the
{\it Riesz Decomposition Property} (RDP) if, for any  $a_{1}$, $a_{2}$, $b_{1}$, $b_{2}\in E$, the equality $a_{1}+ a_{2}=b_{1}+ b_{2}$ implies the existence of four elements $c_{11}, c_{12}, c_{21}, c_{2} \in E$ such that $a_{i}=c_{i1}+ c_{i2}$, and $b_{j}=c_{1j}+ c_{2j}$ for all
$i,j\in \{1,2\}$. }
\end{definition}

We note that due to \cite{DvuPul00}, an effect algebra $(E;+,0,1)$ has  RDP iff, for  $a$, $b_{1}$, $b_{2} \in E$ with $a\leqslant b_{1}+ b_{2},$ there exist  $a_{1},$ $a_{2}\in E$ such that $a =a_{1}+ a_{2}$, and $a_{i}\leqslant b_{i}$ for all $i=1, 2$.

\begin{definition}\label{def:t-norm1}
{ \rm \cite{Good86} An  Abelian po-group $(G;+,0)$ has
the {\it Riesz Decomposition Property} (RDP) if, for any $a$, $b_{1}$, $b_{2} \in G^{+}$ with $a\leqslant b_{1}+ b_{2},$ there exist $a_{1},$ $a_{2}\in G^{+}$ such that $a =a_{1}+ a_{2}$, and
$a_{i}\leqslant b_{i}$ for all $i\in\{1, 2\}$.}
\end{definition}

\section {Effect algebras with RDP and Abelian po-groups}

Let $M$ be a subset of a po-group $G$. We denote by $sss(M)$ the sub-semigroup of $G$  consisting of all finite sums of elements $M$ and of $0.$ An element $u\in G^+$ is said to be  (i) a {\it strong unit} or an {\it order unit} if, given $g \in G,$ there is an integer $n \ge 1$ such that $g  \leqslant nu,$ and (ii) a {\it generative unit} if $G^+= sss(G^+[0,u])$ and $G=G^+- G^-.$ By \cite[Lem 1.4.6]{DvuPul00}, every generative unit is an order unit.

In this section, we give sufficient and necessary conditions such that a po-group $G$ with a generative unit satisfies  RDP.

\begin{definition}\label{de:sequence}
{ \rm Let  $E$ be an atomic effect algebra and $A(E)$ be the set of
atoms of $E$.

{\rm (i)}
Two finite sequences of atoms in $A(E)$ $(a_{1},\ldots,
a_{n})$ and $(b_{1},\ldots,b_{n})$ are called {\it
similar} if there exists a permutation $(p_{1},\ldots, p_{n})$
of $(1,\ldots, n)$ such that $a_{i}=b_{p_{i}},$ $i=1, \ldots, n$.

{\rm (ii)} We say that   $E$ fulfils   the {\it unique atom
representable property} (UARP, for short) if, for any two finite  sequences of atoms
$(a_{1}, \ldots, a_{m} )$ and $(b_{1}, \ldots,
b_{n} )$ such that $\sum_{i=1}^{m}a_{i}=\sum_{j=1}^{n}b_{j},$
then $m=n$ and the sequences
$(a_{1},\ldots,a_{n})$ and $(b_{1},\ldots,b_{n})$ are similar.
}
\end{definition}

Similarly, we can give the following definition for  Abelian po-groups.

\begin{definition}\label{de:sequence1}
{ \rm Let $G$ be   an Abelian po-group and let $A(G^{+})$ be the set of atoms of $G^{+}.$

{\rm (i)} Two finite sequences $(a_{1},\ldots,
a_{n})$ and $(b_{1},\ldots,b_{n})$ of atoms in $A(G^+)$
are called {\it
similar} if there exists a permutation $(p_{1},\ldots, p_{n})$
of $(1,\ldots, n)$ such that $a_{i}=b_{p_{i}},$ $i=1, \ldots, n$.

{\rm (ii)} We say that   $G$ fulfils the {\it
unique atom representable property} (UARP, for short) if, for any   two finite  sequence of atoms
$(a_{1}, \ldots, a_{m})$ and $(b_1,\ldots,b_n)$ in  $A(G^{+}) $  such that $\sum_{i=1}^{m}a_{i}= \sum_{j=1}^{n}b_{j},$ then $m=n$ and the sequences
$(a_{1},\ldots,a_{m})$ and $(b_{1},\ldots,b_{m})$ are similar. }
\end{definition}

\begin{proposition}\label{le:atoms}
Let $G$ be an Abelian po-group with a fixed element $u>0.$
If $G^{+}[0,u]$ satisfies the  condition
$G^{+}[0,u]+G^{+}[0,u]= G^{+}[0,2u]$, then $A(
G^{+}[0,2u])=A(G^{+}[0,u])$, where $A(G^{+}[0, u])$ and
$A(G^{+}[0,2u])$ refer to  the sets of  atoms of the effect algebras
$G^{+}[0, u]$ and $G^{+}[0,2u]$, respectively.
\end{proposition}

\begin{proof} Assume  $a\in A(G^{+}[0,u])$, $b\in G^{+}[0,2u]$ and $b< a$.  Then $b <a \leqslant u,$ so that $b\in G^+[0,u]$ and $a =0.$

Conversely, assume that $a\in A(G^{+}[0,2u])$. Then there exist two
elements $b,$ $c\in  G^{+}[0,u] $ such that $a=b+ c,$  and so $b,
c\leqslant a.$ Since $a\in A(G^{+}[0,2u])$, we have that either $b=0$ or
$c=0$, and so $a\in G^{+}[0,u],$ which implies that $a\in
A(G^{+}[0,u]).$
\end{proof}

\begin{proposition}\label{pr:similar}
Let  $E$ be an effect algebra with RDP. Let  $A=(a_{1},\ldots,a_{n})$ and $B=(b_{1}, \ldots, b_{m})$ be two  finite sequences of atoms such that
$\sum_{i=1}^{n}a_{i}=\sum_{j=1}^{m}b_{j}$, then  the
sequences $A$ and  $B$ are similar.
\end{proposition}

\begin{proof} If $\sum_{i=1}^{n}a_{i}=\sum_{j=1}^{m}b_{j},$ by \cite[Lem 3.9]{DvVe01a}, there is a system $\{x_{ij}\mid i=1,\ldots, n,\ j=1,\ldots, m\}$ of elements from $E$ such that
$$ a_i =\sum_{j=1}^m x_{ij},\quad b_j=\sum_{i=1}^n x_{ij}
$$
for each $i=1,\ldots,n$ and each $j=1,\ldots,m.$ Therefore, for any atom $a_i$ there is a unique $x_{ij_i}$ such that $a_i=x_{ij_i}$ and for any atom $b_j$ there is a unique $x_{i_jj}$ such that $b_j=x_{i_jj}.$ Hence, $n=m$ and the commutativity of $+$ entails $A$ and $B$ are similar.
\end{proof}

\begin{proposition}\label{pr:representation}
Let $G$ be a po-group with RDP and $u$ be
a generative unit, and  let $E=G^{+}[0,u]$ be
an atomic effect algebra. If, for any $x\in E,$ there exists a finite sequence of atoms  $a_{1},\ldots, a_{n}$ in $E$ such that
$x=a_{1}+\cdots+ a_{n}$, then the po-group
$G$ fulfils UARP.
\end{proposition}

\begin{proof} Firstly, the   set $A(G^{+})=\{a\mid a$ is atom of
$G^{+}$$\}$ equals the set $A(E)=\{a\mid a$ is atom of $E\}.$ Since
$E$ is atomic, we have $A(E)\neq\emptyset.$  For any $a\in A(E),$ if
$b\in G^{+}$ with $b< a,$ then we have that $b<u,$ which implies
that $b=0,$ and so $a\in A(G^{+}).$ Conversely, if $a\in A(G^{+}),$
then $a\in G^{+},$ which implies that there exist $a_{1},$ $\ldots,$
$a_{n}\in E$ such that $a=a_{1}+\cdots+a_{n}$. Since $a$ is an atom
of $G^{+}$, we have that there exists a unique index $i\in \{1, \ldots, n\}$
such that $a=a_{i}$ and $a_{j}=0$ with $j\neq i.$ Hence, $a\in E,$
thus $a\in A(E).$

By $G^{+}=ssg(E)$, for any $g\in G^{+},$   there exist   $e_{1},
e_{2}, \ldots, e_{s}\in E$ such that $g=e_{1}+e_{2}+ \cdots+e_{s}$.
Furthermore, by the assumptions, for any $i\in \{1,2,\ldots,s\}$, there
exists a finite sequence of atoms $a_{i1}, a_{i2}, \ldots,  a_{it_{i}}\in E$
such that $e_{i}=a_{i1}+ a_{i2}+ \cdots+ a_{it_{i}},$ and so there
exists a finite sequence of atoms $a_{11}, a_{12}, \ldots,  a_{1t_{1}},\ldots,
a_{s1}, a_{s2}, \ldots,$ $a_{st_{s}}\in E$ such that
$g=a_{11}+a_{12}+ \cdots+a_{1t_{1}}+\cdots+ a_{s1}+ a_{s2}+ \cdots+
a_{st_{s}}$. The rest part of the result follows the similar proof
of Proposition \ref{pr:similar}.
\end{proof}

\begin{theorem}\label{pr:double}
Let $G$ be a po-group $G$ fulfilling UARP
and $u$ be a generative unit.  Then the following
statements hold.

\begin{itemize}
\item[{\rm (i)}] $G^{+}[0,u]$ satisfies RDP.

\vspace{-2mm}\item[{\rm (ii)}] For any natural $n\geqslant 1,$ the effect algebra
$G^{+}[0, nu]$ satisfies RDP.

\vspace{-2mm}\item[{\rm (iii)}] $G^{+}[0,nu]= \underbrace{G^{+}[0,u]+\cdots +
G^{+}[0,u]}_{n-times}$.

\vspace{-2mm}\item[{\rm (iv)}]  The po-group $G$ satisfies RDP.
\end{itemize}
\end{theorem}

\begin{proof}

 (i)\ Assume that  $x\leqslant y+ z$ for  $x,y,z\in G^{+}[0,u].$
 Then there exists an element  $w\in G^{+}[0,u]$ such that
$x+ w=y+ z.$ Since $G$ satisfies UARP, there exist finite sequences of atoms $(x_{1}, \ldots, x_{m}),$ $(w_{1}, \ldots,
w_{q}), $  $(y_{1}, \ldots, y_{n})$ and $(z_{1}, \ldots, z_{p})$
such that $x=x_{1}+ \cdots+ x_{m},$ $w=w_{1}+
\cdots+ w_{q},$ $y=y_{1}+ \cdots+ y_{n}$ and
$z=z_{1}+ \cdots+ z_{p},$ and so  $ x_{1}+
\cdots+ x_{m}+ w_{1}+ \cdots+ w_{q}=y_{1}+
\cdots+ y_{m}+ z_{1}+ \cdots+ z_{p}.$ Hence, the
sequences $(x_{1}, \ldots, x_{m}, w_{1}, \ldots, w_{q})$ and
$(y_{1}, \ldots, y_{n}, z_{1}, \ldots, z_{p})$ are similar, thus for
any $i\in \{1,2,\ldots,m\}$ there exists a unique  $y_{p(i)}$ or a
unique $z_{q(i)}$ such that $x_{i}=y_{p(i)}$ or $z_{q(i)}$. Set
$I_{1}=\{i|$ there exists $y_{p(i)}$ such that $x_{i}=y_{p(i)}\}$,
$I_{2}=\{i|$ there exists $z_{q(i)}$ such that $x_{i}=z_{q(i)}\}$
and we get $a=\sum_{i\in I_{1} }y_{p(i)}$, $b=\sum_{i\in
I_{2}\backslash I_{1} }z_{q(i)}$. Thus, we have that $x=a+ b$
and $a\leqslant y$, $b\leqslant z.$

(ii)\ In any rate,
$G^{+}[0,u]\subseteq G^{+}[0,nu],$ and so $G^{+}=ssg(G^{+}[0,nu])$ which yields that also $nu$ is a generative unit.
By (i), we have that the effect algebra $G^{+}[0,nu]$ satisfies
RDP.

(iii)\ It is easy to see that $G^{+}[0,nu]\supseteq
\underbrace{G^{+}[0,u]+\cdots + G^{+}[0,u]}_{n-times}$.  By (ii),
the effect algebra $G^{+}[0,nu]$ satisfies RDP, and so, for any
$x\in G^{+}[0,nu],$ there exist $n$ elements $x_{1},x_{2},
\ldots,x_{n}$ such that $x=x_{1}+ x_{2}+ \cdots+
x_{n},$ which implies $x\in \underbrace{G^{+}[0,u]+\cdots +
G^{+}[0,u]}_{n-times}$.

(iv)\ For any $a,b,c,d\in G^{+},$  if $a+b=c+d$, then there exists a
natural number $n$ such that $a+b\leqslant nu$. By (ii),
$G^{+}[0,nu]$ satisfies RDP, which implies there exist $x_{1},$
$x_{2},$  $x_{3},$ $x_{4}\in G^{+}[0,nu],$ such that $a=
x_{1}+x_{2},$  $b= x_{3}+x_{4},$ $c= x_{1}+x_{3},$  $d=
x_{2}+x_{4}.$
\end{proof}

An easy corollary of Theorem \ref{pr:double} is the following result.

\begin{corollary}\label{co:uarp}
Let $G$ be a po-group with a generative unit
$u$ and let $E=G^+[0,u]$  If, for any $x\in E,$ there exists a finite sequence of atoms $a_{1},\ldots, a_{n}$ in $E$ such that
$x=a_{1}+\cdots+ a_{n}$.  Then  $G$ satisfies RDP if and only if $G$
satisfies UARP.
\end{corollary}

We say that a poset $E$ satisfies the {\it Riesz Interpolation Property}
(RIP), or $G$ is with {\it interpolation}, if $a_1,a_2 \leqslant b_1,b_2,$ then there is an element $c \in E$ such that $a_1,a_2\leqslant c \leqslant b_1,b_2.$ Then an Abelian po-group $G$ satisfies RIP iff $G$ satisfies RDP, iff $G^+$ satisfies the same property as an effect algebra with RDP,  see \cite[Prop 2.1]{Good86}.

\begin{example}\label{ex:countable}
{\rm \cite{DvuPul00} Let  $G$ be the Abelian group $\mathbb{Z}^{2}$
with the positive cone $G^{+}=\{(a,b)\in G|2a\geqslant b\geqslant 0\}$.

(i) $G$ does not fulfill RIP.

Set $x_{1}=(0,0)$ and $x_{2}=(0,1)$, while $y_{1}=(1,1)$ and
$y_{2}=(1,2)$. Then $x_{i}\leqslant y_{j}$ for all $i,j,$ but there
is no element $z\in G$ such that $x_{i}\leqslant z\leqslant y_{j}$
for all $i,j.$

(ii) The element $u=(2,1)$ is a strong unit of $G$.

 For any
$(a,b)\in G,$ there exists positive element $m$ such that
$(a,b)\leqslant m(2,1)=(2m ,m)$. Notice that $(a,b)\leqslant
n(2,1)=(2n,n)$ iff   $  4n-2a\geqslant n-b\geqslant 0$ iff
$3n\geqslant 2a-b,$ $ n\geqslant b.$ Let $n_{0}=\max\{1, b,
[\frac{1}{3}(2a-b)]+1\}$. Now, we set $m=n_{0}$, then we have that
$m\geqslant 1,$ and so $2m\geqslant m= \max\{1, b,
[\frac{1}{3}(2a-b)]+1\},$ hence,   the inequality $(a,b)\leqslant
m(2,1)=(2m ,m)$ holds.

For any  $(a,b),(c,d)\in G,$ there exist two positive
integers  $m_{1},m_{2}$ such that $(a,b)\leqslant
m_{1}(2,1),(c,d)\leqslant m_{2}(2,1).$ Set $m=\max\{m_{1},m_{2}\},$
we have that $(a,b),(c,d)\leqslant m (2,1).$ Hence, we have prove
that $G$ is directed and the positive element $(2,1)$ is a strong
unit.

(iii) By (ii), the po-group $G$ is directed and we have
$G=G^{+}-G^{+}$.

(iv) Let $0$ and $u$ denote the elements $(0,0)$ and $(2,1),$
respectively. Then the set $G^{+}[0,u]=\{0,(1,0),(1,1),u\}$ is an
interval effect algebra satisfying RDP.

(v) Observe that $G^{+} =\bigcup_{n\in \mathbb{N}}G^{+}[0,nu]$ and
$G^{+} \neq ssg(G^{+}[0,u])$, and so $u$ is not a generative strong
unit   for $G$. The po-group $G$ is not an ambient group with  order
unit $u$ for $E$ (we say that a po-group $G$ with a generative unit $u$ is {\it ambient} for an effect algebra $E$ if $E$ is isomorphic to $G^+[0,u]$).

$G^{+}[0,2u]=\{0,(1,0),(2,0),(3,0),(1,1),(2,1),(3,1),(1,2),(2,2),(3,2),(4,2)\}$,

$G^{+}[0,2u]\neq G^{+}[0,u]+G^{+}[0,u]$. Notice that $(1,2)\in
G^{+}[0,2u]\subseteq G^{+},$ but  for any natural number $n\geqslant 1$, there
exist no elements $x_{i}\in G^{+}[0,u],$ $ i=1,\ldots, n,$ such that
$(1,2)=x_{1}+\cdots+x_{n}$, that is  $(1,2)\notin
ssg(G^{+}[0,u])$.

(vi) Although  $G^{+}[0,u]$ is a Boolean algebra, the effect algebra
$G^{+}[0,2u]$ does not satisfy neither RDP nor  RIP.

For example, $(3,0)+(1,2)=(3,1)+(1,1),$ however, there do not exist any
elements $x_{1},x_{2},x_{3},x_{4} \in G^{+}[0,2u]$ such that
$(3,0)=x_{1}+x_{2},(1,2)=x_{3}+x_{4}$ and
$(3,1)=x_{1}+x_{3},(1,1)=x_{2}+x_{4}$.

For $(2,0),(2,1)\leqslant(3,1),(3,2),$ there exists no element
$x\in G^{+}[0,2u]$ such that  $(2,0), (2,1)\leqslant
x\leqslant(3,1),(3,2).$ }
\end{example}

\begin{example}\label{ex:rdp2}
{\rm Let $G$ be the Abelian group $\mathbb{Z}$, and $G^{+}$ be the
set $\{n\in \mathbb{Z}\mid n=0,$ or $n\geqslant 2\}$. Then $G^{+}$
is a strict cone, and so $G$ is   a po-group with the
partially order $\leqslant_{1}$, for any $a, b\in G$,
$a\leqslant_{1} b$ iff $b-a \in G^{+}$. Let $u=5$, then it is
easy to see that the positive element $u$ is a strong  unit of $G$
and $G^{+}=ssg(E)$, where $E=G^{+}[0,5].$

The equation $G^{+}[0,nu]=\underbrace{G^{+}[0,u]+\cdots
+G^{+}[0,u]}_{n-times}$ holds for any natural number $n$.

The interval effect algebra $G^{+}[0,5]$ is isomorphic to the Boolean
algebra $2^{2}$ which satisfies RDP. $G^{+}[0,10]$ does not
satisfies RDP. In fact, $3+3=2+4,$ however, there exist no
elements $x_{1},$ $x_{2},$ $x_{3},$ $x_{4}\in G^{+}[0,10]$ such that
$3=x_{1}+ x_{2},$ $3=x_{3}+ x_{4},$ $2=x_{1}+ x_{3},$
$4=x_{2}+ x_{4}.$

For any natural number $n\geqslant 2,$ the effect algebra $G^{+}[0,5n]$ does
not fulfil RIP. In fact, $3,4,6,7\in G^{+}[0,5n],$ with $3, 4\leqslant_{1} 6, 7$, however, there is no element $i\in G^{+}[0,5n]$ such that  $3, 4\leqslant_{1}i\leqslant_{1} 6, 7.$
 }
\end{example}

\begin{remark}\label{re:eaRDP}
{\rm  (i)\ Let $G$ be a po-group with the positive cone $G^{+}$, and let a positive element $u$ be a strong unit.
Then the equation $G^{+}=ssg(G^{+}[0, u])$  does not hold in general. See
 Example {\ref{ex:countable}}.

(ii)\ Let $G$ be a po-group with the positive cone $G^{+}$ and let a positive element $u$ be a strong  unit. Then the equation $G^{+}[0,nu]=\underbrace{G^{+}[0, u]+\cdots+G^{+}[0,
u]}_{n-times}$  does not hold, in general. See the   Example {\ref{ex:countable}}.

(iii)\  Let $G$ be a po-group with the positive
cone $G^{+}$, and a positive element $u$ be a strong unit. Assume that
the positive cone $G^{+}=ssg(G^{+}[0, u])$  and
$G^{+}[0,nu]=\underbrace{G^{+}[0, u]+\cdots+G^{+}[0, u]}_{n-times}$
for any natural number $n\geqslant1$. However,  Example {\ref{ex:rdp2}}
shows that  although the effect algebra $G^{+}[0, u]$ satisfy  the
RDP, the effect algebra $G^{+}[0, 2u]$   does not fulfils RDP, which
implies that po-group $G$ does not fulfils  RDP. }
\end{remark}

\section  {Orthocomplete atomic  effect algebra with RDP }

In the present section, we show that every orthocomplete atomic effect algebra with RDP is an MV-effect algebra.

We recall that two elements $a$ and $b$ of an effect algebra $E$ are {\it compatible}, if there exist three elements $a_1,b_1,c\in E$ such that $a= a_1 +c,$ $b=b_1 +c$ and $a_1+b_1+c$ is defined in $E.$ We say that a lattice-ordered effect algebra $E$ is an {\it MV-effect algebra} if all elements of $E$ are mutually compatible.

It is known that if a lattice-ordered effect algebra $E$ satisfies
RDP, then it is also an MV-effect algebra \cite{Rie00}. 

Now, we prove that chain finite effect algebras with RDP are
MV-algebras. Firstly, we recall some useful results for effect
algebras with RDP.

\begin{lemma}\label{le:finite}
{\rm\cite{BeFo95}} Let $E$ be an  effect algebra with RDP. If $E$
is a finite set, then  $E$ is an MV-effect algebra.
\end{lemma}

\begin{definition}\label{de:finite}
{ {\rm\cite{DvuPul00}} Let $E$ be an  effect algebra. If every chain
in $E$ is a finite set, then we say that   $E$ satisfies
the  {\it chain condition.}}
\end{definition}

\begin{lemma}\label{le:fchain}
{\rm\cite{FouBen94}} If an effect algebra $E$  satisfies the chain
condition, then every nonzero element in $E$ is a finite orthogonal
sum of atoms.
\end{lemma}

\begin{theorem}\label{pr:chain}
If an effect algebra $E$ with RDP satisfies the chain condition, then

\begin{itemize}
\vspace{-2mm}\item[{\rm (i)}] $E$ is a finite set.

\vspace{-2mm}\item[{\rm (ii)}] $E$ is an MV-effect algebra.
\end{itemize}
\end{theorem}

\begin{proof}
By Lemma  \ref{le:finite}, it suffices to prove that the statement
(i) holds. Since $E$ satisfies the chain condition,  then there
exists a finite  sequence of atoms  $A=(x_{1},x_{2},\ldots,
x_{n})$ such that $1=x_{1}+ x_{2}+ \cdots+ x_{n}.$
By Proposition \ref{pr:similar}, for any other sequence of atoms
$B=(b_{1},b_2,\ldots, b_{m})$ such that $1=b_{1}+
b_{2}+ \cdots+ b_{m},$ we have that these two
sequences of atoms $A$ and $B$ are similar. Now, for any atom $a$, we have
that $a+ a^{\prime}=1.$ There exists a sequence of atoms $C=(c_{1},c_{2}, \ldots,c_{m})$ such that $a^{\prime}=c_{1}+ c_{2}+ \cdots+ c_{m},$
which implies that the sequence $(a,c_{1},c_{2},\ldots,c_{m})$ is similar to the sequence $A=(x_{1},x_{2},\ldots,x_{n})$. Hence, $a\in \{x_{i}\mid i=1,\ldots,n\}$. Thus, the set of atoms
of $E$ equals $\{x_{i}\mid i=1,\ldots,n\}$. Therefore, for
any $ x\in E$ with $x\neq0$, $x\leqslant 1=x_{1}+ x_{2}+
\cdots+ x_{n}.$ By RDP, there exists a finite sequence of atoms $y_{1},$
$y_{2},$ $\ldots,$ $y_{m}$ with $m\leqslant n$ such that
$x=y_{1}+ y_{2}+ \cdots+ y_{m},$ where
$y_{j}\in\{x_{i}\mid i=1,\ldots,n\}$ for any $j=1, 2, \ldots, m.$
Hence, there exists at most $2^{n}$ elements in $E,$ which implies
that $E$ is a finite set.
\end{proof}

In general,  effect algebras with the chain condition but without RDP are not necessarily finite as the following example  shows.

\begin{example}\label{ex:infi}
{\rm Assume that  $E$ is the horizontal sum of a system $(E_{i})_{i\in \mathbb{N}}$  of effect algebras,  where $E_{i}=\{0,a_{i},1\}$ is a three-element chain effect algebra for any $i\in \mathbb{N}.$ Then $E$ is a chain finite  atomic effect algebra without RDP, and it is also infinite.}
\end{example}

In the following, we  prove that
atomic $\sigma$-orthocomplete effect algebras with RDP are also
MV-effect algebras.

Let $E$ be an effect algebra. We say that a finite sequence $F:=(a_{1},a_{2},\ldots, a_{n})$
is  {\it  orthogonal} if $a_{1}+ a_{2}+ \cdots +
a_{n}$ exists in $E$, and then we write $a_{1}+ a_{2}+
\cdots + a_{n}=\sum_{i=1}^{n}a_{i},$ and the element
$\sum_{i=1}^{n}a_{i}$ is called  the {\it  sum} of the finite system
$F$. The sum of the system $F$ is written as $\sum F.$

For an arbitrary system $A=(a_{i})_{i\in I}$ of not necessarily
different elements of $E$, we say that $G$  is  {\it  orthogonal}, if
every finite subsystem $F$ of $A$ is orthogonal. Furthermore, for an
arbitrary orthogonal system $A$, if the supremum $\bigvee\{\sum
F\mid F$ is a finite subsystem of $ A\}$ exists in $E$, then we say
that the element $\bigvee\{\sum F\mid F$ is a   finite subsystem
of $A\}$ is the {\it  sum} of $A.$ The sum of the system $A$ is
written as $\sum A.$

We say that an effect algebra is {\it  orthocomplete} if an arbitrary
orthogonal  system has a sum. Especially, we say that an effect
algebra is $\sigma$-{\it  orthocomplete} if every countable
orthogonal system  has a sum.

\begin{remark}\label{de:orthogonalcomplete}
 {\rm By \cite[Thm 3.2]{JenPul03},  an effect algebra $E$ is
$\sigma$-orthocomplete iff, for any countable increasing chain
$(a_{i})_{i\in \mathbb{N}}$, the supremum $\bigvee_{i\in \mathbb{N}}
a_{i}$ exists in $E.$}
\end{remark}

\begin{theorem}\label{pr:atompro}
 Let $E$ be a  $\sigma$-orthocomplete atomic effect algebra with
RDP and $A(E)=\{a_{i}\mid i\in \mathbb{N}\}$ be the set of all atoms of
$E$. Then the following statements hold.

\begin{itemize}
\vspace{-2mm}\item[{\rm (i)}]   For any $a_{i}, a_{j}\in A(E)$ with $a_{i}\neq a_{j}$,
then $a_{i}+ a_{j}$ and $ a_{i}\vee a_{j}$ exists and
$a_{i}+ a_{j}=a_{i}\vee a_{j}.$

\vspace{-2mm}\item[{\rm(ii)}]  For any natural number $n\geqslant 2$, the finite set of mutually different atoms $
\{a_{1}, \ldots, a_{n}\}\subseteq A(E)$ is orthogonal in $E$ and
$\sum_{i=1}^{n}a_{i}=\bigvee_{i=1}^{n}a_{i}.$

\vspace{-2mm}\item[{\rm(iii)}]  The system $A(E)$ is an orthogonal system, and
$\sum A(E)= \bigvee A(E)$.
\end{itemize}
\end{theorem}

\begin{proof} (i) See the Lemma 3.2 (ii) in \cite{DvVe01a}.

(ii)  We will proceed by mathematical induction with respect to $n.$

For $n=2$, by (i), $ \{a_{1}, a_{2}\}$ is orthogonal and $
a_{1}+ a_{2}=a_{1}\vee a_{2}.$

Assume that the statement holds for any $m^{\prime}<m.$ For a finite
set of mutually different atoms $ \{a_{1}, \ldots, a_{m}\}$,  by induction hypothesis, we have
that $\sum_{i=1}^{m-1}a_{i}=\bigvee_{i=1}^{m-1}a_{i}.$ Noticing
that for any $i\in \{1, \ldots,m-1\},$ $a_{i}+ a_{m}$ exists,
and so $(\bigvee_{i=1}^{m-1}a_{i})+ a_{m}$ exists, which
implies that the sum $\sum_{i=1}^{m}a_{i}$ exists. Now it
suffices  to prove that
$\sum_{i=1}^{m}a_{i}=\bigvee_{i=1}^{m}a_{i}.$

{\it  Assertion}: If $x\leqslant a_{m}$ and $x\leqslant \bigvee_{i=1}^{m-1}a_{i},$
then $x=0.$

Since $a_{m}$ is an atom and if $x\leqslant a_{m}$, then $x=0$ or
$x=a_{m}.$ Assume that $x=a_{m},$ then there exists an element $b_{1}$ such that
$a_{m}+ b_{1}=(\sum_{i=1}^{m-2}a_{i})+ a_{m-1}.$ By RDP and $a_{m}\wedge a_{m-1}=0$, we will get that
$a_{m}\leqslant \sum_{i=1}^{m-2}a_{i}$. Then there exists an element
$b_{2}$ such that $a_{m}+
b_{2}=(\sum_{i=1}^{m-3}a_{i})+ a_{m-2}.$ Repeating the
same process, we will find  an element $b$ such that
$a_{m}+ b=a_{1}+ a_{2}$, by RDP, we get that
$a_{m}\leqslant a_{1}$ or $a_{n}\leqslant a_{2}$, which is a
contradiction. Consequently, $x=0.$

Now, we assume that $ a_{m},\bigvee_{i=1}^{m-1}a_{i}\leqslant x.$
Then there exist $a, b\in E$ such that $a_{m}+
a=(\bigvee_{i=1}^{m-1}a_{i})+ b=x.$ Using RDP and the
assertion, we have that $a_{m}\leqslant b$, which implies that $
\sum_{i=1}^{m}a_{i}=(\bigvee_{i=1}^{m-1}a_{i})+
a_{m}\leqslant u.$ Hence,
$\sum_{i=1}^{m}a_{i}=\bigvee_{i=1}^{m}a_{i}.$

(iii)\ By (ii), the system $A(E)=\{a_{i}\mid i\in \mathbb{N}\}$ is
an orthogonal system. Since $E$ is a $\sigma$-orthocomplete  effect
algebra, we have that $\sum A(E)= \bigvee\{\sum F\mid F $
is a finite subset of $A(E)\}=\bigvee\{\sum_{i=1}^{n} a_{i}\mid
n\in \mathbb{N},$ $n\geqslant1\}=\bigvee_{i\in \mathbb{N}} a_{i}.$
\end{proof}

The following result generalizes an analogous  result from \cite{DvuCho00}.

\begin{theorem}\label{pr:centrallement}
Let $E$ be a  $\sigma$-orthocomplete atomic effect algebra with
RDP and $A(E)=\{a_{i}\mid i\in I\}$ be the set of all atoms of
$E$ which is at most countable. Let $\imath_{i}$ be the isotropic index of $a_i\in A(E).$ The following statements hold.

\begin{itemize}
\vspace{-2mm}\item[{\rm (i)}]  For any $a_{i}\in A(E)$, the isotropic index
$\imath_{i}$ is finite, $i\in I$.

\vspace{-2mm}\item[{\rm(ii)}] For any $a_{i}\in A(E),$ the interval $E[0,\imath_{i}
a_{i}]=\{x\in E\mid 0\leqslant x\leqslant \imath_{i}a_{i}\}$ equals
to $\{0, a_{i},\ldots, \imath_{i}a_{i}\}.$

\vspace{-2mm}\item[{\rm(iii)}] For any two distinct elements $a_{i}, a_{j}\in A(E),$
$(\imath_{i}a_{i})\wedge(\imath_{j}a_{j})$ exists and
$(\imath_{i}a_{i})\wedge(\imath_{j}a_{j})=0$.

\vspace{-2mm}\item[{\rm(iv)}] For any two distinct elements $a_{i}, a_{j}\in A(E),$
$(\imath_{i}a_{i})+(\imath_{j}a_{j})$ exists and
$(\imath_{i}a_{i})+(\imath_{j}a_{j})=(\imath_{i}a_{i})\vee(\imath_{j}a_{j})$.

\vspace{-2mm}\item[{\rm(v)}] The system $\{\imath_{i}a_{i}\mid a_{i}\in A(E)\}$ is an
orthogonal system, and $\sum\{\imath_{i}a_{i}\mid a_{i}\in
A(E)\}= \bigvee\{\imath_{i}a_{i}\mid a_{i}\in A(E)\}=1$.

\end{itemize}
\end{theorem}

\begin{proof}
(i)\ For any  $a_{i}\in A(E)$, if the sum $na_i$ exists for any
natural number $n\ge 1$, then we get an infinite chain $a_i<2a_i<\cdots<na_i < \cdots <1$.
Since the effect algebra $E$ is $\sigma$-orthocomplete, then
$\bigvee_{n}na_i$ exists. Let $x=\bigvee_{n}na_i$, then $x= \bigvee_{n}(n+1)a_i= a_i+(\bigvee_{n}na_i)= a_i+ x$ which implies that
$a_i=0.$ This is a contradiction with the definition of $a_{i}.$ Hence,
the isotropic $\imath_{i}$ of $a_{i}$ is finite.

(ii)\ For any  $ x\in E[0,\imath_{i} a_{i}]$, if  $x=0$ or
$x=\imath_{i} a_{i}$, then the result holds. Now, if $0<x<\imath_{i}
a_{i}$, then there  exists  $y\in E$ such that $x+ y=\imath_{i}
a_{i}$. By RDP, there exist $x_{11}, \ldots, x_{1i}\in E$ and $x_{21},
\ldots, x_{2i}\in E$ such that $a_{i}=x_{11}+
x_{21}=\cdots=x_{1i}+ x_{2i},$ and $x=x_{11}+\cdots+
x_{1i}$, $y=x_{21}+\cdots+ x_{2i}$. Since $a_{i}$ is an
atom of $E$, we have that $x_{11}, \ldots, x_{1i}, x_{21}, \ldots,
x_{2i}\in \{0, a_{i}\}$, which implies that there exists a natural
number $1\leqslant n\leqslant\imath_{i},$ such that $x=na_{i}.$

(iii)\ For any $x\in E$ with $x\leqslant \imath_{i}a_{i},$
$\imath_{j}a_{j}$, we have that $x\in \{0, a_{i},\ldots,
\imath_{i}a_{i}\}\cap\{0, a_{j},\ldots, \imath_{j}a_{j}\}=0,$ which
implies that $(\imath_{i}a_{i})\wedge(\imath_{j}a_{j})$ exists and
$(\imath_{i}a_{i})\wedge(\imath_{j}a_{j})=0$.

(iv)\ Without loss of generality, we just prove that
$(\imath_{1}a_{1})+(\imath_{2}a_{2})$ exists and
$(\imath_{1}a_{1})+(\imath_{2}a_{2})=(\imath_{1}a_{1})\vee(\imath_{2}a_{2})$.

Noticing that
$(\imath_{1}a_{1})+(\imath_{1}a_{1})^{\prime}=
(\imath_{2}a_{2})+(\imath_{2}a_{2})^{\prime},$
by RDP, there exist $x_{1},$ $x_{2},$ $x_{3},$ $x_{4}\in E$ such
that $\imath_{1}a_{1}=x_{1}+ x_{2},$
$(\imath_{1}a_{1})^{\prime}=x_{3}+ x_{4},$
$\imath_{2}a_{2}=x_{1}+ x_{3},$
$(\imath_{2}a_{2})^{\prime}=x_{2}+ x_{4},$  which implies
$x_{1}=0$ by (iii). Hence, $\imath_{1}a_{1}=x_{2}\leqslant
(\imath_{2}a_{2})^{\prime},$ which implies that
$(\imath_{1}a_{1})+ (\imath_{2}a_{2})$ exists in $E$.
Furthermore, assume that  $\imath_{1}a_{1},$
$\imath_{2}a_{2}\leqslant u,$ for  $u\in E$. Then there exist
$u_{1},$ $u_{2}\in E$ such that $(\imath_{1}a_{1})+
u_{1}=(\imath_{2}a_{2})+ u_{2}$, again by (iii) and RDP, we get
that $\imath_{1}a_{1}\leqslant u_{2}$, which implies that
$(\imath_{1}a_{1})+ (\imath_{2}a_{2})\leqslant u.$ Hence, we
have that $(\imath_{1}a_{1})+
(\imath_{2}a_{2})=(\imath_{1}a_{1})\vee (\imath_{2}a_{2}).$

(v)\ By (iv), the system  $\{\imath_{i}a_{i}\mid a_{i}\in A(E)\}$ is
an orthogonal system, and   $\sum\{\imath_{i}a_{i}\mid a_{i}\in
A(E)\}= \bigvee\{\imath_{i}a_{i}\mid a_{i}\in A(E)\}$. Now,
obviously $\sum\{\imath_{i}a_{i}\mid a_{i}\in A(E)\}\leqslant
1.$ If $\sum\{\imath_{i}a_{i}\mid a_{i}\in A(E)\}< 1,$ then
there exists an element $x\in E$ such that
$(\sum\{\imath_{i}a_{i}\mid a_{i}\in A(E)\})+ x=1$. Since
$E$ is atomic, there exists an atom $a_{x}\leqslant x$. Hence, we
have that $(\sum\{\imath_{i}a_{i}\mid a_{i}\in A(E)\})+
a_{x}$ exists in $E$. But $a_{x}\in A(E),$ we have that
$a_{x}+(\imath_{a_{x}}a_{x})$ exists in $E$, which is a
contradiction. Consequently, we have that $\sum\{\imath_{i}a_{i}\mid
a_{i}\in A(E)\}= 1.$
\end{proof}

We recall that  {\it central elements} of an effect algebra were defined in \cite[Def 1.9.11]{DvuPul00}. If $C(E)$ is the set of central elements, then $C(E)$ is a Boolean algebra.  If $E$ satisfies RDP, then an element $e \in E$ is central iff $e\wedge e'=0,$  see \cite[Thm 3.2]{dv03}.

\begin{lemma}\label{le:orthogonalandjion}
{\rm \cite{JenPul03}} Let $E$ be an orthocomplete  effect algebra.
Let $(a_{\alpha}: \alpha\in \Sigma)\subseteq E$ be an orthogonal
family of central elements. Let $(x_{\alpha}: \alpha\in \Sigma)$ be
a family of elements satisfying $x_{\alpha}\leqslant a_{\alpha},$
for all $\alpha\in \Sigma.$ Then $\bigvee(x_{\alpha}: \alpha\in
\Sigma)$ exists and equals  $\sum(x_{\alpha}: \alpha\in
\Sigma).$
\end{lemma}

\begin{lemma}\label{le:product}
{\rm  \cite{JenPul03}} Let $E$ be an orthocomplete effect algebra.
Let $(a_{\alpha}\mid \alpha\in\Sigma)$ be an orthogonal family of
central elements. Denote $a=\sum(a_{\alpha}\mid
\alpha\in\Sigma).$ Then the element $a$ is  central and $E[0,a]$ is
isomorphic to the product $\prod_{\alpha\in\Sigma
}E[0,a_{\alpha}].$
\end{lemma}


\begin{theorem}\label{pr:atomicsharp}
Let $E$ be a  $\sigma$-orthocomplete  atomic effect algebra with
RDP and $A(E)=\{a_{i}\mid i\in I\}$ be the   set of all atoms of
$E$ which at most countable.  Let $\imath_{i}$ be the isotropic index of $a_i\in A(E).$ Then the following statements hold.

\begin{itemize}

\vspace{-2mm}\item[{\rm (i)}] For any $a_{i}\in A(E),$ the element $\imath_{i}a_{i}$ is a central elements.

\vspace{-2mm}\item[{\rm (ii)}]  For any $a_{i}\in A(E),$ the element $\imath_{i}a_{i}$
is an atom of the Boolean algebra $C(E)$.

\vspace{-2mm}\item[{\rm (iii)}]  For any  $y\in E$, $y=\sum\{y\wedge
\imath_{i}a_{i}\mid a_{i}\in A(E)\}.$

\vspace{-2mm}\item[{\rm (iv)}]  The effect algebra $E$ is isomorphic to the product effect algebra
$\prod_{i\in I} E[0,\imath_{i}a_{i}]$.

\vspace{-2mm}\item[{\rm (v)}] The effect algebra $E$ is a $\sigma$-complete MV-effect algebra.
\end{itemize}
\end{theorem}

\begin{proof}

(i)  By \cite[Thm 3.2]{dv03}, it suffices to prove that
$\imath_{i}a_{i}\wedge (\imath_{i}a_{i})^{\prime}=0.$ Assume that
$x\leqslant  \imath_{i}a_{i},$ $(\imath_{i}a_{i})^{\prime}$. If
$x\neq 0,$ then $a_{i}\leqslant(\imath_{i}a_{i})^{\prime}$ by
Theorem \ref{pr:centrallement} (ii), and so
$a_{i}+(\imath_{i}a_{i})$ exists, which  is a contradiction.
Thus, $x=0,$ and so
$\imath_{i}a_{i}\wedge(\imath_{i}a_{i})^{\prime}=0.$

(ii)  For any $x\in E,$ $x<\imath_{i}a_{i}$, we have that $x\in \{0,
a_{i},\ldots, (\imath_{i}-1)a_{i}\}$ by Theorem
\ref{pr:centrallement} (ii). If $x\neq 0,$ then $a_{i}\leqslant x,
x^{\prime}$ by Theorem \ref{pr:centrallement} (ii), which
implies that $x\notin C(E).$ Hence, we have that $\imath_{i}a_{i}$
is an atom of $C(E).$

(iii)  Since $\imath_{i}a_{i}\in C(E)$, for any $y\in E,$ $y\wedge
\imath_{i}a_{i}$ exists in $E$. Since the set $\{\imath_{i}a_{i}\mid
a_{i}\in A(E)\}$ is orthogonal, the set  $\{y\wedge
\imath_{i}a_{i}\mid a_{i}\in A(E)\}$ is also orthogonal and so the
sum  $\sum\{y\wedge \imath_{i}a_{i}\mid a_{i}\in A(E)\}$ exists
in $E$. Notice that for any two elements $y\wedge \imath_{i}a_{i}$,
$y\wedge \imath_{j}a_{j}\in \{y\wedge \imath_{i}a_{i}\mid a_{i}\in
A(E)\},$ the sum $(y\wedge \imath_{i}a_{i})+(y\wedge
\imath_{j}a_{j})$ exists and  $(y\wedge
\imath_{i}a_{i})+(y\wedge \imath_{j}a_{j})=(y\wedge
\imath_{i}a_{i})\vee(y\wedge \imath_{j}a_{j})$ by Lemma
\ref{le:orthogonalandjion}. Whence, for any finite subset $F\subseteq
\{y\wedge \imath_{i}a_{i}\mid a_{i}\in A(E)\},$ we have that
$\sum\{x\mid x\in F\}=\bigvee\{x\mid x\in F\}$. In addition, we have
that $\sum\{y\wedge \imath_{i}a_{i}\mid a_{i}\in A(E)\}=
\bigvee\{y\wedge \imath_{i}a_{i}\mid a_{i}\in A(E)\}\leqslant y.$
Assume that  $\sum\{y\wedge \imath_{i}a_{i}\mid a_{i}\in
A(E)\}= \bigvee\{y\wedge \imath_{i}a_{i}\mid a_{i}\in A(E)\}< y.$
Then there exists an element $x\in E$ such that
$x+(\sum\{y\wedge \imath_{i}a_{i}\mid a_{i}\in A(E)\})=y,$
and so there exists an atom $a_{i_{0}}\in E$ such that
$a_{i_{0}}\leqslant x.$ However, $x+(y\wedge
\imath_{i_{0}}a_{i_{0}})$ exists, and so $a_{i_{0}}+(y\wedge
\imath_{i_{0}}a_{i_{0}})\leqslant y,\imath_{i_{0}}a_{i_{0}},$ hence,
$a_{i_{0}}+(y\wedge \imath_{i_{0}}a_{i_{0}})\leqslant
y\wedge(\imath_{i_{0}}a_{i_{0}}),$ which is a contradiction. Thus,
we have that $\sum\{y\wedge \imath_{i}a_{i}\mid a_{i}\in
A(E)\}= \bigvee\{y\wedge \imath_{i}a_{i}\mid a_{i}\in A(E)\}= y.$

(iv) By Theorem \ref{pr:centrallement} (v) and Lemma
\ref{le:product}, the statement holds.

(v) For any $i\in I$, the chain $E[0,\imath_{i}a_{i}]$ is a
finite MV-effect algebra, and so it is complete.
Hence, the product
$\prod_{i\in I}E[0,\imath_{i}a_{i}]$ is   also a $\sigma$-complete
MV-effect algebra.
\end{proof}

\begin{remark}\label{re:BAD-POSETS}
{\rm  In Proposition 3.12 of \cite{DvuCho00}, the authors proved
that any atomic $\sigma$-complete Boolean D-poset with the
countable set of atoms  $\{a_{i}\mid i\in I\}$ can be expressed
as a direct product of finite chains. In fact, any Boolean D-poset is
an MV-effect algebra, which is also a lattice-ordered effect algebra
with RDP \cite{DvuPul00}. By Theorem \ref{pr:atomicsharp}, we
can see that any $\sigma$-orthocomplete  atomic effect algebra with
RDP is also a lattice-ordered, thus it is an MV-algebra.
Furthermore, similar to Theorem \ref{pr:atomicsharp}, we can
prove the following result.}
\end{remark}

\begin{theorem}\label{pr:atomicorthocomp}
Let $E$ be an  orthocomplete  atomic effect algebra with  RDP and
$A(E)=\{a_{i}\mid i\in I\}$ be the set of  atoms of $E$. Then the
following statements hold.

\begin{itemize}
\vspace{-2mm}\item[{\rm (i)}] For any $a_{i}\in A(E),$ the element $\imath_{i}a_{i}$ is
a central elements.

\vspace{-2mm}\item[{\rm (ii)}]  For any $a_{i}\in A(E),$ the element $\imath_{i}a_{i}$
is an atom of Boolean algebra $C(E)$.

\vspace{-2mm}\item[{\rm (iii)}] For any  $y\in E$, $y=\sum\{y\wedge
\imath_{i}a_{i}\mid a_{i}\in A(E)\}.$

\vspace{-2mm}\item[{\rm (iv)}]  The effect algebra $E$ is isomorphic to the product effect algebra
$\prod_{i\in I} E[0,\imath_{i}a_{i}]$.

\vspace{-2mm}\item[{\rm (v)}] The effect algebra $E$ is a complete MV-effect algebra.
\end{itemize}
\end{theorem}

\section{Applications}

In the present section, we apply the methods and the results of the previous sections to a noncommutative generalization of effect algebras, pseudo-effect  algebras, and to a description of the state space of some effect algebras.

A noncommutative generalization of effect algebras was introduced in \cite{DvVe01a,DvVe01b} and some additional basic properties can be found in \cite{dv03}.

\begin{definition}\label{def:pea}
{\rm  \cite{DvVe01a} A structure
$(E;+,0,1),$  where $+$ is a partial binary operation and 0 and 1 are constants, is called a {\it pseudo-effect  algebra},  if for all $a,b,c \in E,$ the following hold.
\begin{itemize}
\item[{\rm (PE1)}] $ a+ b$ and $(a+ b)+ c $ exist if and only if $b+ c$ and $a+( b+ c) $ exist, and in this case,
$(a+ b)+ c =a +( b+ c)$.

\item[{\rm (PE2)}] There are exactly one  $d\in E $ and exactly one $e\in E$ such
that $a+ d=e + a=1 $.

\item[{\rm (PE3)}] If $ a+ b$ exists, there are elements $d, e\in E$ such that
$a+ b=d+ a=b+ e$.

\item[{\rm (PE4)}] If $ a+ 1$ or $ 1+ a$ exists,  then $a=0$ .
\end{itemize}}
\end{definition}
We recall that a pseudo-effect  algebra $E$ is an effect algebra iff the partial addition $+$ is commutative, i.e. $a+b$ exits in $E$ iff $b+a$ is defined in $E,$  and the $a+b = b+a.$

In the same way as for effect algebras, we define for pseudo-effect  algebra (i) the isotropic index, $\imath(a),$ of any element $a$ of a pseudo-effect  algebra, (ii) atom, (iii) atomic system, (iv) RDP, (v) central element, and (vi) center $C(E).$

We say that a pseudo-effect  algebra $E$ is {\it monotone $\sigma$-complete} provided that every ascending sequence $x_1 \leqslant   x_2 \leqslant \cdots $ of elements  in $E$ has a supremum $x = \bigvee_n x_n.$  We recall that if $E$ is an effect algebra, then the notions $\sigma$-orthocomplete effect algebras and monotone $\sigma$-complete effect algebras are equivalent. For pseudo-effect algebras, the notion of  the $\sigma$-orthocomplete pseudo-effect algebra is not straightforward in view of the non-commutativity of the partial addition $+.$ Therefore, for our aims, we prefer the notion of the monotone $\sigma$-completeness of pseudo-effect algebras.

\begin{theorem}\label{pr:atomicorthocomp1}
Let $E$ be a monotone $\sigma$-complete  atomic pseudo-effect algebra with  RDP and
$A(E)=\{a_{i}\mid i\in I\}$ be the set of  atoms of $E$ that is at most countable. Then $E$ is a commutative PEA, i.e., $E$ is an effect algebra.
\end{theorem}

\begin{proof}  For any two atoms $a,b$ with $a\neq b,$ we have $a\wedge b = 0,$ so that by \cite[Lem 3.1]{dv03}, $a+b,$   $b+a$ and $a\vee b$ exists in $E$ and they are equal.

We assert that the isotropic index of any atom $a$ of $E$ is finite. Indeed, assume the converse, i.e. $\imath(a)=\infty.$  Then $na \in E$ for any integer $n\geqslant 1,$ and $\bigvee_n na \in E.$  Hence, $\bigvee_n(n+1)a = \bigvee_n na + a$ implies the contradiction $a=0.$

In the same way as in Theorem \ref{pr:centrallement}(i) we prove that every $\imath(a)a$ is a central element.

Since $A(E)$ is at most countable, we assume that $A(E)=\{a_1,a_2,\ldots\}.$
The RDP implies that, for all $a,b \in A(E),$ $\imath(a)a \wedge \imath(b)b=0$, which yields  that $\imath(a)a+\imath(b)b=\imath(b)b+\imath(a)a=\imath(a)a\vee\imath(b)b.$  In the same way, we can show that if $a_1,\ldots,a_n$ are mutually different atoms, then $\imath(a_1)a_1+\cdots + \imath(a_n)a_n$ exists in $E$ and it equals $b_n:=\imath(a_1)a_1\vee \cdots\vee \imath(a_n)a_n.$ In addition, $b_n = \imath(a_{j_1})a_{j_1}+\cdots + \imath(a_{j_n})a_{j_n}$ for any permutation $(j_1,\ldots,j_n)$ of $(1,\ldots,n).$
Thus, we have that $\{b_n\}$ is an ascending sequence and so  $\bigvee_n b_n$ exists in $E$ and we claim that $\bigvee_n b_n=1.$ In fact,  if $\bigvee_n b_n<1,$ then there exists an atom $a$ such that $\bigvee_n b_n + a\leqslant 1,$ and so,  $\imath(a)a+a$ exists in $E$ which is absurd. Hence, $\bigvee_n \imath(a_n) a_n = \bigvee_n b_n=1.$ By \cite[Thm 5.11]{dv03}, we have $x = \bigvee_n (x\wedge \imath (a_n) a_n)$ for any $x \in E.$ Therefore, by \cite[Pro 6.1(ii)]{dv03}, there exists an isomorphism $\phi:E\rightarrow \prod_{i\in I}E[0,\imath(a_{i})a_{i}],$ where $\phi(x)=(x\wedge \imath(a_{i})a_{i})_{i\in I}.$ Further, for any $x\in E$ and any $ i\in \mathbb{N},$ we have $x\wedge\imath(a_{i})a_{i}\in E[0,\imath(a_{i})a_{i}]= \{0,a,\ldots,\imath(a_{i})a_{i}\}$ by RDP.

For any $x,y\in E,$   $x+y$ exists in $E$ iff $(x\wedge\imath(a_{i})a_{i})+(y\wedge\imath(a_{i})a_{i})$ exists. Thus, if  $x+y$ exists in $E,$ then we have that for any $i\in I,$
 $(x+y)\wedge\imath(a_{i})a_{i}=$
$(x\wedge\imath(a_{i})a_{i})+(y\wedge\imath(a_{i})a_{i}) =(y\wedge\imath(a_{i})a_{i})+(x\wedge\imath(a_{i})a_{i}) $
which implies that $y+x$ exists and $x+y=y+x.$
\end{proof}

Now we apply Theorem \ref{pr:atomicorthocomp} for the description of states on some atomic effect algebras.

We say that a {\it state} on an effect algebra $E$ is a mapping $s:E \to [0,1]$ such that (i) $s(a+b) = s(a)+s(b)$ whenever $a+b$ is defined in $E$, and (ii) $s(1)=1$. A state is an analogue of a probability measure.  A state $s$ is said to be  {\it extremal} if, for any
states $s_{1},s_{2}$ and $\alpha\in(0,1)$, the equation $s=\alpha s_{1}+(1-\alpha)s_{2}$ implies $s=s_{1}=s_{2}.$ Let $\mathcal S(E)$ and $\partial_e \mathcal S(E)$ denote the set of all states and extremal states, respectively, on $E.$  We recall that it can happen that an effect algebra is stateless. But every interval effect algebra admits at least one state, see \cite[Cor 4.4]{Good86}. We say that a net of states $\{s_\alpha\}_\alpha$ {\it converges weakly} to a state $s$ iff $\lim_\alpha s_\alpha(a)=s(a)$ for any $a \in E.$  Then $\mathcal S(E)$ is a compact Hausdorff space, and due to the Krein--Mil'man Theorem, see e.g. \cite[Thm 5.17]{Good86}, every state is a weak limit of a net of convex combinations of extremal states on $E.$

We recall that a state on an MV-effect algebra is extremal, \cite{Mun}, iff $s(a\wedge b)=\max\{s(a),s(b)\}$ for all $a,b \in E.$

A state $s$ is $\sigma$-additive if for any monotone sequence $\{a_i\}$ such that $\bigvee_i a_i =a$ implies $s(a)=\lim_i s(a_i).$ Equivalently, if $a=\sum_n a_n,$ then $s(a)=\sum_n s(a_n).$

\begin{theorem}\label{th:5.3}
Let $E$ be a  $\sigma$-orthocomplete atomic effect algebra with
RDP and $A(E)=\{a_{i}\mid i\in I\}$ be the set of all atoms of
$E$ that is at most countable.  Let $\imath_{i}$ be the isotropic index of $a_i\in A(E).$ For any $i\in I,$ we define a mapping $s_i:E \to [0,1]$ via

$$ s_i(a) = \max\{j\mid  ja_i \leqslant a  \wedge \imath_{i} a_i\}/ \imath_{i},\quad a \in E.
$$
Then $s_i$ is  an extremal state on $E$ which is also $\sigma$-additive. If $s$ is a $\sigma$-additive state on $E,$ then $s(a) = \sum_i \lambda_i s_i(a),$ $a \in E.$ Moreover, every extremal state that is also $\sigma$-additive is just of the form $s_i$ for a unique $i,$ and a state $s=s_i$  for some $i\in I$ if and only if $s(\imath_ia_i)=1.$
\end{theorem}

\begin{proof}  By Theorem \ref{pr:atomicorthocomp}(i),(iii), the element $\imath_{i} a_i$ is central and $a = \sum_i\{a\wedge \imath_{i}a_i\}.$  Therefore, $s_i(a)$ is a real number from the real interval $[0,1]$ and $(a+b)\wedge \imath_{i}a_i = (a\wedge \imath_{i}a_i)+(b\wedge \imath_{i}a_i)$ which proves that $s_i$ is a state. Since by  Theorem \ref{pr:atomicorthocomp}(v), $E$ is an MV-effect algebra. If $a,b \in E,$ we have $(a\wedge b)\wedge \imath_{i} a_i) = (a\wedge \imath_{i} a_i)\wedge (b \wedge \imath_{i} a_i)$ which implies $s_i(a\wedge b)=\min\{s_i(a),s_i(b)\}$ which proves $s_i$ is an extremal state.

By \cite[Thm 5.11]{dv03}, if $x_n \nearrow x $ and $e$ is a central element, then $(\bigvee_n x_n)\wedge e = \bigvee_n(x_n \wedge e).$  From this and the definition of $s_i,$ we have that each $s_i$ is $\sigma$-additive.

Let $s$ be an arbitrary $\sigma$-additive state, then $a = \sum_i\{a\wedge \imath_{i}a_i\}$ and $1= \sum_i \imath_{i}a_i,$ so that $s(a) = \sum_i s(a\wedge \imath_{i}a_i) = \sum_i \lambda_i s_i(a),$ where $\lambda_i = s(\imath_{i}a_i).$

Therefore, if $s$ an extremal state that is also $\sigma$-additive, from the previous decomposition we conclude that $s=s_i$ for a unique $i.$

Now assume that $s$ is a state on $E$ such that $s(\imath_i a_i)=1.$  Then $s(a_i)=1/\imath_i.$  Since $\imath_i a_i$ is a central element, for any $a \in E,$ we have $s(a) = s(a\wedge \imath_ia_i) + s(a\wedge (\imath_ia_i)')= s(a\wedge \imath_ia_i) = s_i(a).$
\end{proof}

\section{Conclusion}

In the paper, we have studied effect algebras $E$ which are also MV-effect algebras, i.e. every two elements of $E$ are compatible. Since every MV-effect algebra satisfies the Riesz Decomposition Property, in other words, every  two decompositions of the unit element $1$ have a joint refinement, we have concentrated to effect algebras with RDP. We recall that RDP fails for $\mathcal E(H),$ and every effect algebra with RDP is an interval in an interpolation Abelian po-group with strong unit, and every MV-effect algebra is an interval in a lattice ordered group with strong unit.

The main result says, Theorem \ref{pr:centrallement}, that every $\sigma$-orthocomplete atomic effect algebra with RDP and with the countable set of atoms is in fact an MV-effect algebra which is the countable direct product of finite chains.

This results was applied also for pseudo-effect algebras, where it was proved, Theorem \ref{pr:atomicorthocomp1}, that any analogous pseudo-effect algebra has to be commutative. In addition, the studied methods allow also to give a complete characterization of $\sigma$-additive states of our type of effect algebras, Theorem \ref{th:5.3}.




\vspace{3mm}
{\bf  Acknowledgement:} A.D. thanks  for the support by Center of Excellence SAS -~Quantum
Technologies~-,  ERDF OP R\&D Project
meta-QUTE ITMS 26240120022, the grant VEGA No. 2/0059/12 SAV and by
CZ.1.07/2.3.00/20.0051  and  MSM 6198959214.

Y.X. thanks  for the support by SAIA,
n.o. (Slovak Academic Information Agency) and the Ministry of
Education, Science, Research and Sport of the Slovak Republic. This
work is also supported by National Science Foundation of China
(Grant No. 60873119),  and the Fundamental Research Funds for the
Central Universities (Grant No. GK200902047).


\begin{thebibliography}{99}

\vspace{-2mm}\bibitem{BeFo95}  M. K. Bennet and D. Foulis,
\emph{Phi-symmetric effect algebras},
Foundations of Physics {\bf 25} (1995), 1699--1722.

\vspace{-2mm}\bibitem{Cig00} R. Cignoli, I. M. L. D'Ottaviano  and D. Mundici,
\emph{``Algebraic Foundations of Many-Valued Reasoning",}
Trends in Logic, Volume 7, Kluwer Academic Publishers, Dordrecht, 2000.


\bibitem{Cha}
 C.C. Chang,
{\it Algebraic analysis of many valued logics,}
 Transactions of the American Mathematical Society
{\bf 88} (1958),
 467--490.

\vspace{-2mm}\bibitem{dv03} A. Dvure\v{c}enskij,
\emph{Central elements and
 Cantor-Bernstein's theorem for pseudo-effect algebras,}
 Journal of the Australian Mathematical Society {\bf  74} (2003), 121--143.



\vspace{-2mm}\bibitem{DvuCho00}  A. Dvure\v{c}enskij,   F. Chovanec, E. Ryb\'{a}rikov\'{a},
    \emph{D-homomorphisms and atomic
$\sigma$-complete Boolean D-posets,} Soft Computing   {\bf 4} (2000),
9--18.

\vspace{-2mm}\bibitem{DvuPul00}  A. Dvure\v{c}enskij and  S. Pulmannov\'{a}.
\emph{``New Trends in Quantum Structures",}
Kluwer Academic Publishers,
Dordrecht, and Ister Science,  Bratislava, 2000.

\vspace{-2mm}\bibitem{DvVe01a} A. Dvure\v{c}enskij, T. Vetterlein
\emph{Pseudoeffect algebas. I. Basic properties,}
International Journal of Theoretical Physics {\bf  40}
(2001), 685--701.


\vspace{-2mm}\bibitem{DvVe01b} A. Dvure\v censkij, T. Vetterlein,
Pseudoeffect algebras. II. Group representation, {\it International Journal of Theoretical Physics} {\bf 40} (2001), 703--726.

\vspace{-2mm}\bibitem{FouBen94}  D. Foulis and M. K. Bennet,
\emph{Effect algebras and unsharp quantum logics,} Foundations of
Physics {\bf 24} (1994),   1325--1346.

\vspace{-2mm}\bibitem{Good86}  K. R. Goodearl,
\emph{``Partially Ordered Abelian Groups with Interpolation",}
Mathematical Surveys and Monographs No 20, American Mathematical
Society, Providence, Rhode Island, 1986.

\vspace{-2mm}\bibitem{Jen04}  G. Jen\v{c}a,
\emph{Boolean algebras R-generated by MV-effect algebras,} Fuzzy
Sets and Systems {\bf 145} (2004),  279--285.

\vspace{-2mm}\bibitem{JenPul03} G. Jen\v{c}a, S. Pulmannov\'{a},
\emph{Orthocomplete effect algebras}, Proceedings of the American
Mathematical Society {\bf  131} (2003), 2663--2671.

\vspace{-2mm}\bibitem{KopCho94}  F. K\^{o}pka  and F. Chovanec,
\emph{D-posets},
Mathematica Slovaca {\bf  44} (1994),  21--34.


\vspace{-2mm}
\bibitem{Mun}
 D. Mundici,
{\it Averaging the truth-value in \L ukasiewicz logic,} Studia
Logica
 {\bf 55} (1995), 113--127.

\vspace{-2mm}\bibitem{Rav98}  K. Ravindran,
\emph{States on effect algebras that have the $\phi$-sysmmetric
property}, International Journal of Theoretical Physics {\bf 37} (1998),
175--181.

\vspace{-2mm}\bibitem{Rie00}  Z. Rie\v{c}anov\'{a},
\emph{Generalization of blocks for D-lattices and lattice ordered
effect algebras}, International Journal of Theoretical Physics {\bf 39} (2000), 231--237.

\end{thebibliography}
\end{document}